\documentclass{ws-procs9x6} 

\usepackage{graphicx} 

\usepackage{amsfonts}
\usepackage{amssymb}
\usepackage{amsmath}

\usepackage{cite}
\usepackage{epsf}

\def\idos{invariant differential operators}

\def\ido{invariant differential operator}

\def\lra{\longrightarrow}

\def\llr{\longrightarrow}
  \def\tV{{\tilde V}}

 \def\cgc{{\cg^\bbc}}
 \def\np{\vfill\eject}
\def\a{\alpha}
\def\b{\beta}
\def\d{\delta}

\def\vr{\vert}

\def\l{\lambda}

\def\D{{\Delta}}
\def\ca{{\cal A}} \def\cb{{\cal B}} \def\cc{{\cal C}}
\def\cd{{\cal D}} \def\ce{{\cal E}} \def\cf{{\cal F}}
\def\cg{{\cal G}} \def\ch{{\cal H}} 
 \def\ck{{\cal K}} 
\def\cm{{\cal M}} \def\cn{{\cal N}} 
\def\cp{{\cal P}} \def\cq{{\cal Q}}  
\def\ct{{\cal T}}

\def\bbz{\mathbb{Z}}
\def\bbc{\mathbb{C}}
\def\bbr{\mathbb{R}}
\def\bbn{\mathbb{N}}

\def\eqn#1{\begin{equation}\label{#1}}
\def\ee{\end{equation}}

\def\bea{\begin{eqnarray}}
\def\eea{\end{eqnarray}}

\def\eqnn#1{\begin{eqnarray}\label{#1}}
\def\eqna#1{\begin{subequations} \label{#1}
\begin{eqnarray}}

\def\eena{\end{eqnarray}
\end{subequations}}

\def\chh{\chi^+}
\def\Lp{\Lambda^+}

\def\sevha{{\textstyle{7\over2}}}

\def\nn{\nonumber}
\def\ha{{\textstyle{1\over2}}}
\def\d{\delta}
\def\L{\Lambda} \def\l{\lambda}

\def\rf#1#2{(\ref{#1}{#2})}

   \def\nt{\noindent}

\def\veps{\varepsilon}
\def\r{{\rho}}

\input epsf.tex
\newcount\figno
\def\fig#1#2#3{
\par\begingroup\parindent=0pt\leftskip=1cm\rightskip=1cm\parindent=0pt
\baselineskip=11pt \global\advance\figno by 1 
\epsfxsize=#3 \centerline{\epsfbox{#2}} \vskip 12pt
#1\par
\endgroup\par}
\def\figlabel#1{\xdef#1{\the\figno}}
\def\encadremath#1{\vbox{\hrule\hbox{\vrule\kern8pt\vbox{\kern8pt
\hbox{$\displaystyle #1$}\kern8pt} \kern8pt\vrule}\hrule}}

\begin{document}

\title{Invariant Differential Operators for the\\[4pt] Real Exceptional Lie Algebra $F''_4$}

\author{V.K. Dobrev}

 \address{Institute for Nuclear Research and Nuclear Energy,\\ Bulgarian
Academy of Sciences,\\ 72 Tsarigradsko Chaussee,  1784 Sofia,
Bulgaria}

 \begin{abstract}
 In the present paper we continue the project of systematic
construction of invariant differential operators on the example of
the non-compact  exceptional Lie algebra  $F''_4$ which is the split rank one form of the exceptional Lie algebra $F_4$.
  We classify
the reducible Verma modules over $F_4$ which are compatible with this induction.
Thus, we obtain the classification of  the corresponding invariant differential operators.
\end{abstract}

\bigskip

\hfill {\it Contribution to Peter Suranyi Festschrift}

\bodymatter

\section{Introduction}

Invariant differential operators   play very important role in the
description of physical symmetries. One of the most important such symmetry is the conformal symmetry. Incidentally,
our involvement in this area started in the 70s with the Euclidean conformal group ~$SO(N,1)$\cite{DMPPT,DPPT,DP}. It is important that this group is of split rank
one.\footnote{Though not
in our main exposition flow, we should mention that in parallel were developed the study and applications of the Minkowskian conformal group $SO(N,2)$, cf., e.g.,  \cite{confmink}.}

Our  general scheme for constructing  \idos\
operators was given some time ago \cite{Dob}.
After our paper was submitted we learned
that similar questions   were treated in a different approach (by so-called K-induction)\cite{Zhel}. There were treated the split rank one groups $SU(N,1)$ and $SO(N,1$).

In recent papers \cite{Dobinv,Dobparab} we started the systematic explicit
construction of the invariant differential operators. Yet the two remaining split rank one cases, namely, $F''_4$ and $Sp(N,1)$, were not treated until now.\footnote{We are
thankful to Joachim  Hilgert for pointing out this omission.}
Thus, it was important to treat the case $F''_4$ which we do in the present paper.

The first task in  our approach is to make the multiplet classification of
the reducible Verma modules over the algebra in consideration\cite{Dobmul}.
Such classification provides the weights of embeddings between the Verma modules via
the singular vectors, and thus\cite{Dob}, the weights of the invariant differential operators.

We have done the multiplet classification for many real non-compact algebras\cite{Dobk1}.
In the present paper we  focus on the complex exceptional Lie algebra ~$F_4$~ and on its
 split  rank one form algebra  $F''_4$.
 Our scheme requires that we use induction from parabolic subalgebras.

 The importance of the parabolic subgroups comes from the fact that
the representations induced from them generate all (admissible)
irreducible representations of $G$\cite{Lan,Zhea,KnZu}.

 In the split rank one case $F''_4$ there is only one nontrivial parabolic:
 ~$\cp = \cm \oplus \ca \oplus \cn$, where ~$\cm = so(7)$, $\dim_\bbr \cn = 15$.

  We present the multiplet classification of
 the reducible Verma modules over $F_4$ which are compatible with the   parabolic $\cp$ of $F''_4$.
 We give also the weights of the singular vectors between
 these modules. By our scheme\cite{Dob} these  singular vectors will produce
 the invariant differential operators.

\section{Preliminaries}

\subsection{Lie algebra}

We start with the complex exceptional Lie algebra  ~$\cg^\bbc ~=~ F_4$. We use
the standard definition of ~$\cg^\bbc$~ given in terms of the
Chevalley generators $X^\pm_i ~, ~H_i ~, ~i=1,2,3,4 (=$rank$\,F_4)$, by the relations~:
\eqnn{com}
 &[H_j\,, \,H_k] \, = \, 0 \,, \,\,\,[H_j\,, \,X^\pm_k] \, = \, \pm
a_{jk} X^\pm_k \,,
\,\,\,[X^+_j \,, \,X^-_k] \, = \,
\d_{jk} \,H_j \,,  \\
&\sum_{m=0}^n \,(-1)^m \,\left({n \atop m}\right)
\,\left(X^\pm_j\right)^m \,X^\pm_k \,\left(X^\pm_j\right)^{n-m}
\,=\, 0 \,, \,\,j \neq k \,, \,\,n = 1 - a_{jk} \,, \nn\eea where
\eqn{aijf4} (a_{ij}) = \begin{pmatrix} 2 & -1 & 0 & 0 \cr -1 & 2
& -1 & 0\cr 0 & -2 & 2 & -1\cr 0 &0 &-1 & 2\end{pmatrix} ~; \ee
is the Cartan matrix of $\cg^\bbc$, ~$\a^\vee_j
\,\equiv\, {2 \a_j\over (\a_j , \a_j)}$ is the co-root of
$\a_j\,$, \,\, $(\cdot , \cdot)$ is the scalar product of the roots,
so that the nonzero products between the simple roots are:
$(\a_1, \a_1) = (\a_2, \a_2) = 2(\a_3, \a_3) = 2(\a_4, \a_4)~=~2 $, ~
$(\a_1, \a_2) = -1$, ~$(\a_2, \a_3) = -1$, ~$(\a_3, \a_4) = -1/2$.
 The
elements $H_i$ span the Cartan subalgebra $\ch$ of $\cg^\bbc$, while
the elements $X^\pm_i$ generate the subalgebras $\cg^\pm$. We shall
use the standard triangular decomposition \eqn{deca} \cg^\bbc \,=\,
\cg_+\oplus\ch\oplus \cg_- \,, \qquad\cg_\pm \,\equiv
\,\mathop{\oplus}\limits_{\a\in\D^\pm} \,\cg_\a \,, \ee where
$\D^+$, $\D^-$, are the sets of positive, negative, roots, resp.
Explicitly    we have that there are roots of two lengths with
length ratio $2:1$.\\ \indent The long roots are:
~$\a_1$, ~$\a_2$,
~$\a_1 + \a_2$,
~$\a_2 + 2\a_3$,
~$\a_1 + \a_2 + 2\a_3$,
~$\a_1 + 2\a_2 + 2\a_3$,
~$\a_2 +2\a_3 + 2\a_4$,
~$\a_1 +\a_2 +2\a_3 +2\a_4$,
~$\a_1 +2\a_2 +2\a_3 +2\a_4$,
~$\a_1 +2\a_2 + 4\a_3 +2\a_4$,
~$\a_1 + 3\a_2 + 4\a_3 + 2\a_4$,
~$2\a_1 + 3\a_2 + 4\a_3 + 2\a_4$.
With the chosen normalization they have length 2.\\
 \indent The short roots are:
~$\a_3$, ~$\a_4$,
~$\a_2 + \a_3$,
~$\a_3 + \a_4$,
~$\a_1 + \a_2 + \a_3$,
~$\a_2 + \a_3 + \a_4$,
~$\a_1 +\a_2 +\a_3 +\a_4$,
~$\a_2 +2\a_3 +\a_4$,
~$\a_1 +2\a_2 +2\a_3+\a_4$,
~$\a_1 +\a_2 +2\a_3 +\a_4$,
~$\a_1 +2\a_2 + 3\a_3 + \a_4$,
~$\a_1 +2\a_2 +3\a_3 +2\a_4$,
and they have length 1.\\
(Note that the short roots are exactly those which contain $\a_3$
and/or $\a_4$ with coefficient 1, while the long roots contain $\a_3$
and $\a_4$ with even coefficients.)

Thus, as well-known, $F_4$ is 52--dimensional ($52 = \vert \D\vert +$ rank $F_4$).

In terms of the normalized basis ~$\veps_1, \veps_2, \veps_3,
\veps_4$ ~we have:
\eqnn{posf4} \D^+ = \{\veps_i, ~1 \leq i \leq 4; ~\veps_j \pm
\veps_k, ~1 \leq j <k \leq 4; \nn\\ {1 \over 2}(\veps_1 \pm \veps_2 \pm
\veps_3\pm \veps_4), {\rm ~all\ signs}\} ~. \eea The simple roots are: \eqn{pif4}
\pi = \{\a_1 = \veps_2 - \veps_3, ~\a_2 = \veps_3 - \veps_4, \a_3 = \veps_4,
~\a_4 = {1 \over 2} (\veps_1 - \veps_2 - \veps_3 - \veps_4)\}
~.\ee

 For later use we note that the Lie algebra ~$B_4=so(9)^\bbc$~ may be embedded most easily in ~$F_4$~
as the Lie algebra generated by the
  16 roots on the first line of \eqref{posf4}. Indeed, the latter form the positive root system
of ~$B_4$~ with simple roots ~$\veps_i-\veps_{i+1}\,, ~i=1,2,3, ~\veps_4\,$.

The Weyl group of ~$F_4$~ is the semidirect product of ~$S_3$~ with
a group which itself is the semidirect product of ~$S_4$~ with
~$(\bbz/2\bbz)^3$, thus, ~$|W| = 2^7\,3^2 = 1152$.

\subsection{Verma modules}

Let us recall that a ~{\it Verma module} ~$V^\L$~ is defined as
the HWM over ~$\cgc$~ with highest weight ~$\L \in \ch^*$~ and
highest weight vector ~$v_0 \in V^\L$, induced from the
one-dimensional representation ~$V_0 \cong \bbc v_0$~ of
~$U(\cb)$~, where ~$\cb  = \ch \oplus \cg_+$~ is a Borel
subalgebra of ~$\cgc$, such that:
\eqnn{indb}
& &X ~v_0 ~~=~~ 0 , \quad \forall\, X\in \cg_+ \cr
&&H ~v_0 ~~=~~ \L(H)~v_0\,, \quad \forall\, H \in \ch \eea

Verma modules are generically irreducible. A Verma  module ~$V^\L$~ is
reducible\cite{BGG} iff there exists a root ~$\b \in\D^+$~ and ~$m\in\bbn$~
such that
\eqn{red} (\L + \r~, ~\b^\vee) ~=~  m \ee
holds, where ~$\r = {1 \over 2}\sum_{\a \in \D^+}~\a$~,
($ \r ~=~ 8\a_1 +15\a_2 +21\a_3 + 11\a_4 $).

If \eqref{red} holds then the reducible Verma module $V^\L$  contains an invariant submodule
which is also a Verma module ~$V^{\L'}$~ with shifted weight ~$\L'=\L-m\b$.
This statement is equivalent to the fact that $V^\L$ contains a
~{\it singular vector}~
~$v_s \in V^\L$, such that ~$v_s ~\neq ~\xi v_0\,$, ($0\neq\xi\in\bbc$),
and~:
\eqnn{sing}
&& X ~v_s ~~=~~ 0 , \quad \forall\, X\in \cg_+ \cr
&&H ~v_s ~~=~~ \L'(H) ~v_s\,, \quad
\L' ~=~ \L - m\b, ~~\forall\,
H \in \ch \eea
More explicitly\cite{Dob},
\eqn{singp} v^s_{m,\b} = \cp_{m,\b}\, v_0 \ . \ee

   The general reducibility conditions \eqref{red} for $V^\L$
   spelled out for the  simple roots in our situation are:
\eqnn{reda}
&& m_1 ~\equiv~ m_{\a_1} ~=~ (\L +\r, \a_1), ~~~
m_2 ~\equiv~  m_{\a_2} ~=~ (\L +\r, \a_2),\\
&& m_3 ~\equiv~  m_{\a_3} ~=~ (\L +\r, 2\a_3), ~~~
m_4 ~\equiv~ m_{\a_4} ~=~ (\L +\r, 2\a_4) \nn\eea
If we write
\eqn{lamm} \L = \l_1\a_1 + \l_2\a_2 + \l_3\a_3 + \l_4\a_4\ee
 then the numbers $\l_i$ are expressed through $m_i$ as follows:
 \eqnn{lammm} &&\l_1 = 2m_1 +3m_2 +2m_3 +m_4 -8 , \nn\\  && \l_2 = 3m_1 +6m_2 +4m_3 +2m_4 -15 , \nn\\ && \l_3 = 24 -5m_1 -10m_2 -6m_3 -3m_4  , \nn\\
 && \l_4 = \ha (23 -5m_1) -5m_2 -3m_3 -m_4 \eea

The numbers ~$m_i$~ from \eqref{reda} corresponding to the simple roots are called Dynkin labels, while the more general Harish-Chandra parameters  are:
\eqn{hcpar} m_\b ~=~ (\L +\r, \b^\vee), ~~~\b\in\D^+ \ee
Note that the expression from \eqref{sing} ~$\L' = \L -m\b$ may be written using \eqref{hcpar} as a Weyl reflection:
\eqn{weylx} s_\b (\L+\r)  ~\doteq~ \L+\r - m_\b\b ~=~ \L'+\r \ee

Explicitly, the Harish-Chandra parameters (for the non-simple roots) in terms of the Dynkin labels $m_1,m_2,m_3,m_4$ are:
\eqnn{hclong}
&& m_{\a_1 + \a_2} = m_1+m_2 ~\equiv ~ m_{12}, \nn\\ &&
m_{\a_2 + 2\a_3} = m_2+m_3 ~\equiv ~ m_{23}, \nn\\ &&
m_{\a_1 + \a_2 + 2\a_3} = m_1+m_2 + m_3~\equiv ~ m_{13} ,\nn\\ &&
m_{\a_1 + 2\a_2 + 2\a_3} = m_1+2m_2 + m_3~\equiv ~ m_{13,2} ,\nn\\ &&
m_{\a_2 +2\a_3 + 2\a_4} = m_2 + m_3 +m_4 ~\equiv ~ m_{24}, \nn\\ &&
m_{\a_1 +\a_2 +2\a_3 +2\a_4} = m_1+m_2 + m_3+ m_4 ~\equiv ~ m_{14} ,\nn\\ &&
m_{\a_1 +2\a_2 +2\a_3 +2\a_4} = m_1+2m_2 + m_3 +m_4 ~\equiv ~ m_{14,2},\nn\\ &&
m_{\a_1 +2\a_2 + 4\a_3 +2\a_4} = m_1+2m_2 +2 m_3+ m_4 ~\equiv ~ m_{14,23},\nn\\ &&
m_{\a_1 + 3\a_2 + 4\a_3 + 2\a_4} = m_1+3m_2 +2 m_3+ m_4 ~\equiv ~ m_{14,23,2},\nn\\ &&
m_{2\a_1 + 3\a_2 + 4\a_3 + 2\a_4}  = 2m_1+3m_2 +2 m_3+ m_4 ~\equiv ~ m_{14,13,2}
\eea
for the long roots, while for the short roots we have:
\eqnn{hcshort}
&& m_{\a_2 + \a_3} = 2m_2+m_3 ~\equiv ~ m_{23,2}, \nn\\ &&
 m_{\a_1 + \a_2 + \a_3} = 2m_1+2m_2+m_3 ~\equiv ~ m_{13,12} , \nn\\ &&
m_{\a_3 + \a_4}  = m_3 +m_4 ~\equiv ~ m_{34}, \nn\\ &&
m_{\a_2 + \a_3 + \a_4} = 2m_2 +  m_3 +m_4 ~\equiv ~ m_{24,2}, \nn\\ &&
m_{\a_1 +\a_2 +\a_3 +\a_4}  = 2m_1+2m_2 + m_3 +m_4 ~\equiv ~ m_{14,12} , \nn\\ &&
m_{\a_2 +2\a_3 +\a_4} = 2m_2 + 2m_3 +m_4 ~\equiv ~ m_{24,23} ,\nn\\ &&
m_{\a_1 +2\a_2 +2\a_3+\a_4} = 2m_1+4m_2 + 2m_3 +m_4 ~\equiv ~ m_{14,13,2,2} ,\nn\\ &&
m_{\a_1 +\a_2 +2\a_3 +\a_4} = 2m_1+2m_2 + 2m_3 +m_4 ~\equiv ~ m_{14,13} ,\nn\\ &&
m_{\a_1 +2\a_2 + 3\a_3 + \a_4} = 2m_1+4m_2 + 3m_3 +m_4 ~\equiv ~ m_{14,13,23,2} ,\nn\\ &&
m_{\a_1 +2\a_2 +3\a_3 +2\a_4} = 2m_1+4m_2 + 3m_3 +2m_4 ~\equiv ~ m_{14,14,23,2}
\eea
where we have introduce short=hand notation, e.g, ~$m_{13,2} \equiv m_{13} +m_2$,
 ~$m_{13,12} \equiv m_{13} +m_{12}$, etc.

\subsection{Structure theory of the real form}

The split real form of ~$F_4$~ is denoted as ~$F''_4\,$,
sometimes as  ~$F_{4(-20)}\,$. It has rank four.
Its maximal compact subalgebra is ~$\ck \cong so(9)$, also of rank four.
This real form has   discrete series
representations since ~rank$F''_4\,~=$ ~rank$\,\ck$.
The number of discrete series is equal to the ratio ~$\vert W(\cgc,\ch^\bbc)\vert / \vert W(\ck^\bbc,\ch^\bbc)\vert$,
where $\ch$ is a compact Cartan subalgebra of both $\cg$ and $\ck$, ~$W$ are the relevant Weyl groups \cite{Knapp}.
Thus, the number of discrete series in our setting is three. They will be identified below.


The Iwasawa decomposition of ~$\cg \equiv F''_4\,$,~ is: \eqn{iwa32}
\cg ~=~ \ck \oplus \ca \oplus \cn \ , \ee the Cartan decomposition
is: \eqn{car32} \cg ~=~ \ck \oplus \cq , \ee where we use
~$\dim_\bbr\,\cq = 16$,
 ~$\dim_\bbr\,\ca = 1$, ~$\cn =  \cn^+\,$, or
 ~$\cn =  \cn^- \cong \cn^+\,$, ~$\dim_\bbr\,\cn^\pm =15$.

   Since ~$\cg$~ is  of split rank one,  thus, the minimal (also maximal) parabolic ~$\cp$~
and the corresponding Bruhat decomposition are:
\eqnn{min32} && \cp ~=~ \cm \oplus \ca \oplus \cn \ , \qquad  \cm = so(7) \\
&& \cg ~=~ \cm \oplus  \ca \oplus \cn^+\oplus \cn^- \nn\eea

Note that the root system of ~$\cm^\bbc  = so(7,\bbc) = B_3$~ consists of the roots
  \eqn{posb3} \D_3^+ = \{\veps_i, ~2 \leq i \leq 4; ~\veps_j \pm
\veps_k, ~2 \leq j <k \leq 4 \}\ee
which are part of \eqref{posf4}, while the
 simple roots are part of \eqref{pif4}
\eqn{poss3} \pi_3 ~=~ \{ \a_1 = \veps_2 -\veps_3, \quad \a_2 = \veps_3 -\veps_4,  \quad \a_3 = \veps_4 \} \ee

The roots of ~$\cm^\bbc$~ are called ~{\it $\cm$-compact roots} of the ~$F_4$~ root system \eqref{posf4},
the rest are called ~{\it $\cm$-noncompact roots}. The latter give rise to \idos, as explained below.

More explicitly, the $\cm$-compact roots are:
\eqna{compr}  && \a_1, ~\a_2,
~\a_1 + \a_2 \equiv \a_{12},
~\a_2 + 2\a_3 \equiv \a_{23,3},
~\a_1 + \a_2 + 2\a_3 \equiv \a_{13,3}  ,\nn\\ &&
~\a_1 + 2\a_2 + 2\a_3\equiv \a_{13,23},  \\
&& \a_3,
~\a_2 + \a_3\equiv \a_{23},
~\a_1 + \a_2 + \a_3\equiv \a_{13} ,
\eena
\rf{compr}a are long roots, \rf{compr}b - short.\\
The $\cm$-noncompact roots are:
\eqna{noncomp}
&& \a_2 +2\a_3 + 2\a_4 \equiv \a_{24,23},
~\a_1 +\a_2 +2\a_3 +2\a_4 \equiv \a_{14,34}, \nn\\ &&
\a_1 +2\a_2 +2\a_3 +2\a_4 \equiv \a_{14,24},\nn\\
&&\a_1 +2\a_2 + 4\a_3 +2\a_4 \equiv \a_{14,24,3,3},
~\a_1 + 3\a_2 + 4\a_3 + 2\a_4\equiv \a_{14,24,23,3} ,\nn\\
&& 2\a_1 + 3\a_2 + 4\a_3 + 2\a_4\equiv \a_{14,14,23,3} \\
&&\a_4, ~\a_3 + \a_4\equiv \a_{34},
~\a_2 + \a_3 + \a_4 \equiv \a_{24},
~\a_1 +\a_2 +\a_3 +\a_4\equiv \a_{14}, \nn\\ &&
\a_2 +2\a_3 +\a_4\equiv \a_{24,3}, ~
\a_1 +2\a_2 +2\a_3+\a_4\equiv \a_{14,23}, \nn\\ &&
~\a_1 +\a_2 +2\a_3 +\a_4 \equiv \a_{14,3},
~\a_1 +2\a_2 + 3\a_3 + \a_4 \equiv \a_{14,23,3}, \nn\\ &&
\a_1 +2\a_2 +3\a_3 +2\a_4\equiv \a_{14,24,3}
 \eena
 \rf{noncomp}a are long roots, \rf{noncomp}b - short.

Correspondingly, the Dynkin labels ~$m_1,m_2,m_3$~ are called $\cm$-compact, while ~$m_4$~ is called
$\cm$-noncompact.

\subsection{Elementary representations}

Further, let ~$G,K,P,M,A,N$~ are Lie groups with Lie algebras  ~$\cg_0,\ck,\cp,\cm,\ca,\cn$.

Let ~$\nu$~ be a (non-unitary) character of ~$A$, ~$\nu\in\ca^*$.
 Let ~ $\mu$ ~ fix a finite-dimensional unitary representation
~$D^\mu$~ of $M$ on the space ~$V_\mu\,$.

 We call the induced
representation ~$\chi =$ Ind$^G_{P}(\mu\otimes\nu \otimes 1)$~ an
~{\it \it elementary representation} of $G$\cite{DMPPT}. (These are
called also {\it generalized principal series representations} (or {\it
limits thereof})\cite{Knapp}.)   Their spaces of functions are:  \eqn{func}
\cc_\chi ~=~ \{ \cf \in C^\infty(G,V_\mu) ~ \vr ~ \cf (gman) ~=~
e^{-\nu(H)} \cdot D^\mu(m^{-1})\, \cf (g) \} \ee where ~$a=
\exp(H)\in A$, ~$H\in\ca\,$, ~$m\in M$, ~$n\in N$. The
representation action is the {\it left regular action}:  \eqn{lrega}
(\ct^\chi(g)\cf) (g') ~=~ \cf (g^{-1}g') ~, \quad g,g'\in G\ .\ee

An important ingredient in our considerations are the
highest/lowest weight representations~ of ~$\cg$. These can be
realized as (factor-modules of) Verma modules ~$V^\L$~ over
~$\cg$, where ~$\L\in (\ch)^*$,   the weight ~$\L = \L(\chi)$~ being determined
uniquely from $\chi$\cite{Dob}.

As we have seen when a Verma module is reducible and \eqref{red} holds then there is a
singular vector \eqref{singp}. Relatedly, then
there exists\cite{Dob} an {\it \ido}  \eqn{invop}  \cd_{m,\b} ~:~ \cc_{\chi(\L)}
~\llr ~ \cc_{\chi(\L-m\b)} \ee given explicitly by: \eqn{singvv}
 \cd_{m,\b} ~=~ \cp_{m,\b}(\widehat{\cn^-})  \ee where
~$\widehat{\cn^-}$~ denotes the {\it right action} on the functions
~$\cf$.

Actually, since our ERs  are induced from finite-dimensional
representations of ~$\cm$~  the corresponding Verma modules are
always reducible. Thus, it is more convenient to use ~{\it
generalised Verma modules} ~$\tV^\L$~ such that the role of the
highest/lowest weight vector $v_0$ is taken by the
(finite-dimensional) space ~$V_\mu\,v_0\,$.

Algebraically, the above is governed by the notion of ~$\cm$-compact roots of $\cgc$.   The consequence
of this is that \eqref{red} is always fulfilled for the ~$\cm$-compact roots of $\cgc$. That is why we consider
generalised Verma modules. Relatedly, the \idos\ corresponding to  ~$\cm$-compact roots are trivial.

 \section{Main multiplets of $F''_4$}

Further  we classify the generalized Verma modules (GVM) relative to the    parabolic subalgebra $\cp$ \eqref{min32}.
This also provides the classification of the $P$-induced ERs  with the same Casimirs.
 The classification is done as follows.   We group the
reducible Verma modules  (also the corresponding ERs) related by nontrivial embeddings
in sets called ~{\it multiplets}\cite{Dobmul,Dob}. These multiplets
may be depicted as a connected graph, the
vertices of which correspond to the GVMs  and the lines
between the vertices correspond to the GVM embeddings (and also the \idos\ between the ERs).
The explicit parametrization of the multiplets and of their Verma modules (and ERs) is
important for understanding of the situation.

The result of our classification is a follows. The multiplets of GVMs (and ERs) induced from ~$\cp$~
 are parametrized by four positive integers - the Dynkin labels.
   Each
multiplet contains  24 GVMs (ERs).
 Their signatures are given as follows:
\eqnn{mult4} && \chi^-_0 ~=~ \{ m_1,m_2,m_3,m_4\} \\
&&\chi^-_a ~=~ \{ m_1,m_2,m_{34},-m_4\}\nn\\
&&\chi^-_b ~=~ \{ m_1,m_{23},m_{4},-m_{34}\}\nn\\ &&
\chi^-_c ~=~ \{ m_{12},m_{23},m_{4},-m_{24,2}\}\nn\\ &&
\chi^-_d ~=~ \{ m_{2},m_{13},m_{4},-m_{14,12}\}  \nn\\ &&
\chi^-_e ~=~ \{ m_{13}, m_2, m_{34}, -m_{24,23}\} \nn\\ &&
\chi^-_f ~=~ \{ m_{23},m_{12},m_{34},-m_{14,13}\}\nn\\ &&
\chi^-_g ~=~ \{ m_{14}, m_2, m_{3}, -m_{24,23}\}\nn\\ &&
\chi^-_h ~=~ \{m_{23},m_{1},m_{24,2},-m_{14,13,2,2}  \}\nn\\ &&
\chi^-_i ~=~ \{ m_{24},m_{12},m_{3},-m_{14,13}\}\nn\\ &&
\chi^-_j ~=~ \{m_{2},m_{1},m_{24,23},-m_{14,13,23,2}  \} \nn\\ &&
\chi^-_k ~=~ \{m_{24},m_{1},m_{23,2},-m_{14,13,2,2}  \}\nn\\ &&
 \chh_k  ~=~ \{ m_{24},m_{1},m_{23,2},-m_{14,13,23,2}  \}  \nn\\ && 
 \chh_j  ~=~ \{m_{2},m_{1},m_{24,23},-m_{14,14,23,2}  \}\nn\\ && 
 \chh_i  ~=~ \{ m_{24},m_{12},m_{3},-m_{14,13,23,2}  \}\nn\\ && 
 \chh_h  ~=~ \{ m_{23},m_{1},m_{24,2},-m_{14,14,23,2}  \} \nn\\ && 
 \chh_g   ~=~ \{ m_{14},m_{2},m_{3},-m_{14,13,23,2}  \} \nn\\ && 
 \chh_f ~=~ \{ m_{23},m_{12},m_{34},-m_{14,14,23,2}  \} \nn\\ && 
  \chh_e  ~=~ \{ m_{13},m_{2},m_{34},-m_{14,14,23,2}  \}\nn\\ && 
 \chh_d  ~=~ \{ m_{2},m_{13},m_{4},-m_{14,14,23,2} \}  \nn\\ && 
\chh_c   ~=~ \{ m_{12},m_{23},m_{4},-m_{14,14,23,2}  \}    \nn\\ && 
\chh_b  ~=~ \{ m_{1},m_{23},m_{4},-m_{14,14,23,2}  \}\nn\\ && 
  \chh_a ~=~ \{ m_{1},m_{2},m_{34},-m_{14,14,23,2}  \}   \nn\\ &&  
 \chh_0 ~=~ \{ m_{1},m_{2},m_{3},-m_{14,14,24,2}  \}  \nn 
  \eea

These multiplets  are presented in  Fig. 1.
On the figure each arrow represents an embedding between two Verma modules, ~$V^\L$ and ~$V^{\L'}$,~ the arrow pointing to the embedded module ~$V^{\L'}$.
Each arrow carries a number ~$n$, $n=1,2,3,4$, which indicates the level of the embedding,  ~$\L' = \L - m_n\,\b$\cite{Dobk1}. By our construction it also
represents the invariant differential operator ~$\cd_{n,\b}\,$, cf. \eqref{invop}.

Further, we note that there is an additional symmetry  w.r.t. to the dashed line in  Fig. 1.
It is relevant for the ERs and indicates the integral intertwining Knapp-Stein (KS) operators\cite{KnSt} acting between the
  spaces ~$\cc_{\chi^\mp}$~ in opposite directions:
\eqn{ksks} G^+_{KS} ~:~ \cc_{\chi^-} \llr \cc_{\chi^+}\ , \qquad
G^-_{KS} ~:~ \cc_{\chi^+} \llr \cc_{\chi^-} \ee
Note that the KS opposites are induced from the same irreps of ~$\cm$.

This symmetry may be more explicit if we change the parametrization:
\eqn{chanp}  \{ m_1,m_2,m_3,m_4\} \lra [ m_1,m_2,m_3; c ] \ee
so that the action of the KS operators on this signature is:
\eqn{ksksi} G^\pm_{KS} ~:~ [ m_1,m_2,m_3; c ] \lra [ m_1,m_2,m_3; -c ]\ee
This enables us to write the multiplet in a more compact way:
\eqnn{mult4c} && \chi^\pm_0 ~=~ [ m_1,m_2,m_3;  \pm (m_{14,2,4} + m_3/2)] \\
&&\chi^\pm_a ~=~ [ m_1,m_2,m_{34};\pm \ha m_{14,13,23,2} ]\nn\\
&&\chi^\pm_b ~=~ [ m_1,m_{23},m_{4}; \pm \ha m_{14,13,2,2} ]\nn\\ &&
\chi^\pm_c ~=~ [ m_{12},m_{23},m_{4}; \pm \ha m_{14,13} ]\nn\\ &&
\chi^\pm_d ~=~ [ m_{2},m_{13},m_{4}; \pm \ha m_{24,23} ]  \nn\\ &&
\chi^\pm_e ~=~ [ m_{13}, m_2, m_{34};  \pm \ha m_{14,12}] \nn\\ &&
\chi^\pm_f ~=~ [ m_{23},m_{12},m_{34}; \pm \ha m_{24,2} ]\nn\\ &&
\chi^\pm_g ~=~ [ m_{14}, m_2, m_{3}; \pm \ha  m_{13,12}  ]\nn\\ &&
\chi^\pm_h ~=~ [m_{23},m_{1},m_{24,2}; \pm \ha  m_{34} ]\nn\\ &&
\chi^\pm_i ~=~ [ m_{24},m_{12},m_{3}; \pm \ha  m_{23,2} ]\nn\\ &&
\chi^\pm_j ~=~ [m_{2},m_{1},m_{24,23}; \pm \ha m_{4}  ] \nn\\ &&
\chi^\pm_k ~=~ [m_{24},m_{1},m_{23,2}; \pm \ha m_{3}  ]\nn
    \eea

Note that if in \eqref{mult4c} we denote generically    \eqn{genc} \chi^\pm = \{ m_1,m_2,m_3,   m_4^\pm \} = [m_1,m_2,m_3 ;  c^\pm] \ee
 then there is the relation
\eqn{genad} |c^+|+|c^-|  = |m^+_4 | + |m^-_4| \ . \ee

\bigskip

Standardly  we take ~$\L^-_0$ as the top Verma module (ER), since it is not embedded in any other Verma module.
Relatedly, the VM $\Lp_0$  has no embedded Verma modules - it is the KS opposite of $\L^-_0$
We have  indicated also the other KS opposites by denoting their signatures with "$\pm$".

{\it Remark:} ~Note that the pairs ~$\chi^\pm_j$~ and ~$\chi^\pm_k$~  are related by KS operators, but in each case the operator ~$G^+_{KS}$~   is degenerated into  a differential operator, namely, we have
  \eqna{multopjk} && \L^-_j ~~ {m_4\a_{14,24,3} \atop \longrightarrow } ~~\Lp_j \\ && \phantom{\L^-_j} \nn\\
&& \L^-_k ~~ {m_3\a_{14,24,3} \atop \longrightarrow } ~~\Lp_k \eena

\section{Reduced multiplets}

\subsection{Main reduced multiplets}

The above multiplets do not exhaust the relevant \idos.  There are reduced multiplets which are obtained as we go to the walls of the relevant Weyl chambers.
First, there are four main reduced multiplets $M_k$, $k=1,2,3,4$, which may be obtained by setting the parameter $m_k=0$.


The main reduced multiplet $M_1$ contains 18   GVMs (ERs).
 Their signatures are given as follows:

\eqnn{mult4r} && \chi^-_0 ~=~ \{ 0,m_2,m_3,m_4\} \\
&&\chi^-_a ~=~ \{ 0,m_2,m_{34},-m_4\}\nn\\
&&\chi^-_b ~=~ \{ 0,m_{23},m_{4},-m_{34}\}\nn\\ &&
\chi^-_c ~=~ \{ m_{2},m_{23},m_{4},-m_{24,2}\} = \chi^-_d \nn\\ &&
 \chi^-_e ~=~ \{ m_{23}, m_2, m_{34}, -m_{24,23}\} = \chi^-_f \nn\\ &&
  \chi^-_g ~=~ \{ m_{24}, m_2, m_{3}, -m_{24,23}\} = \chi^-_i\nn\\ &&
\chi^-_h ~=~ \{m_{23},0,m_{24,2},-m_{24,23,2,2}  \}\nn\\ &&
 \chi^-_j ~=~ \{m_{2},0,m_{24,23},-m_{24,23,23,2}  \} \nn\\ &&
\chi^-_k ~=~ \{m_{24},0,m_{23,2},-m_{24,23,2,2}  \}\nn\\ &&
 \chh_k  ~=~ \{ m_{24},0,m_{23,2},-m_{24,23,23,2}  \} \nn\\ && 
  \chh_j    ~=~ \{m_{2},0,m_{24,23},-m_{24,24,23,2}  \} \nn\\ && 
 \chh_h  ~=~ \{ m_{23},0,m_{24,2},-m_{24,24,23,2}  \}  \nn\\ &&  
  \chh_g  ~=~ \{ m_{24},m_{2},m_{3},-m_{24,23,23,2}  \} = \chh_i \nn\\ &&  
\chh_e    ~=~ \{ m_{23},m_{2},m_{34},-m_{24,24,23,2}  \}= \chh_f  \nn\\ && 
\chh_c   ~=~ \{ m_{2},m_{23},m_{4},-m_{24,24,23,2}  \} =   \chh_d\nn\\ &&  
\chh_b  ~=~ \{ 0,m_{23},m_{4},-m_{24,24,23,2}  \} \nn\\ &&  
\chh_a   ~=~ \{ 0,m_{2},m_{34},-m_{24,24,23,2}  \}    \nn\\ && 
 \chh_0  ~=~ \{ 0,m_{2},m_{3},-m_{24,24,24,2}  \}  \nn 
  \eea

  Note that only six of these GVMs (ERs), namelty, ~$\chi^-_c,\chi^-_e,\chi^-_g, \chh_c,\chh_e,\chh_g$, are induced  by finite-dimensional representations
  of $\cm$. We give the latter also in the more compact notation:

  \eqnn{mult4c1}
&& \chi^\pm_c ~=~ [ m_{2},m_{23},m_{4}; \pm \ha m_{24,23} ]\\ &&
\chi^\pm_e ~=~ [ m_{23}, m_2, m_{34};  \pm \ha m_{24,2}] \nn\\ &&
 \chi^\pm_g ~=~ [ m_{24}, m_2, m_{3}; \pm \ha  m_{23,2}  ]\nn
      \eea

    These GVMs are related in the following way:
  \eqn{multop41} \begin{matrix}  \L^-_c & {m_3\a_{24} \atop \longrightarrow } & \L^-_e &  {m_4\a_{24,3} \atop \longrightarrow } & \L^-_g \cr
&&&&\cr \updownarrow  && \updownarrow && \updownarrow \cr
 &&&&\cr    \Lp_c & {m_3\a_{14,3} \atop \longleftarrow } & \Lp_e &  {m_4\a_{14} \atop \longleftarrow } & \Lp_g
\end{matrix}
\ee
where the up-down arrows designate the $G^\pm$ KS operators.

The reduced multiplets of type $M_2$ also contain 18 GVMs (ERs). We give only the signatures of those that are induced  by finite-dimensional representations  of $\cm$~:

\eqnn{mult42}
 &&\chi^-_c ~=~ \{ m_1,m_{3},m_{4},-m_{34}\} = \chi^-_b \nn\\ &&
 \chi^-_f ~=~ \{ m_{3},m_{1},m_{34},-m_{1,34,1,3}\} =\chi^-_h \nn\\ &&
 \chi^-_k ~=~ \{ m_{34},m_{1},m_{3},-m_{1,34,1,3}\} =\chi^-_i  \nn\\ &&
   \chh_k  ~=~ \{ m_{34},m_{1},m_{3},-m_{1,34,1,3,3}  \}=   \chh_i \nn\\ &&
 \chh_f  ~=~ \{ m_{3},m_{1},m_{34},-m_{1,34,1,34,3}  \}=   \chh_h\nn\\ &&
    \chh_c  ~=~ \{ m_{1},m_{3},m_{4},-m_{1,34,1,34,3,2}  \} =  \chh_b\nn
   \eea
or in the more compact notation:
\eqnn{mult4c2}  &&
\chi^\pm_c ~=~ [ m_{1},m_{3},m_{4}; \pm  ( m_{13} + m_4/2 ) ]\\ &&
\chi^\pm_f ~=~ [ m_{3},m_{1},m_{34}; \pm \ha m_{34} ]\nn\\ &&
 \chi^\pm_k ~=~ [ m_{34},m_{1},m_{3}; \pm \ha  m_{3} ]\nn
    \eea

These GVMs are related in the following way:
  \eqn{multop42} \begin{matrix}  \L^-_c & \phantom{\lra}  & \L^-_f &  {m_4\a_{24,3} \atop \longrightarrow } & \L^-_k \cr
&&&&\cr \updownarrow  && \updownarrow && \updownarrow \cr
 &&&&\cr    \Lp_c &  & \Lp_f &  {m_4\a_{4} \atop \longleftarrow } & \Lp_k
\end{matrix}
\ee
Note that the Remark before \eqref{multopjk} is again valid for the pair $\chi^\pm_k$.

The reduced multiplets of type $M_3$  contain 15 GVMs (ERs).  Those induced  by finite-dimensional representations
  of $\cm$~ are:
\eqnn{mult43}
&&\chi^-_b ~=~ \{ m_1,m_2,m_{4},-m_4\} = \chi^-_a \\ &&
 \chi^-_c ~=~ \{ m_{12},m_{2},m_{4},-2m_{2}-m_4 \} = \chi^-_e \nn\\ &&
\chi^-_d ~=~ \{ m_{2},m_{12},m_{4},-m_{12} -m_4 \}= \chi^-_f \nn\\ &&
 \chi^-_j ~=~ \{m_{2},m_{1},m_{2,2,4},-m_{12,12,2,2,4}  \} = \chi^-_h \nn\\ &&
\chi^-_k ~=~ \{m_{2,4},m_{1},m_{2,2},-m_{12,12,2,2,4}  \}= \chh_k \nn\\ &&
  \chh_j  ~=~ \{m_{2},m_{1},m_{2,2,4},-m_{12,12,2,2,4,4}  \}=  \chh_h\nn\\ && 
 \chh_d  ~=~ \{ m_{2},m_{12},m_{4},-m_{12,12,2,2,4,4}  \}= \chh_f  \nn\\ && 
  \chh_c ~=~ \{ m_{12},m_{2},m_{4},-m_{12,12,2,2,4,4}  \}= \chh_e \nn\\ &&  
 \chh_b  ~=~ \{ m_{1},m_{2},m_{4},-m_{12,12,2,2,4,4}  \}  = \chh_a   \nn  
    \eea
or
\eqnn{mult4c3}
&&\chi^\pm_b ~=~ [ m_1,m_{2},m_{4}; \pm  (m_1 + 2m_2 + m_4/2)    ]\\ &&
\chi^\pm_c ~=~ [ m_{12},m_{2},m_{4}; \pm  ( m_{12} + m_4/2)  ]\nn\\ &&
\chi^\pm_d ~=~ [ m_{2},m_{12},m_{4}; \pm  ( m_{2} + m_4/2)  ]  \nn\\ &&
 \chi^\pm_j ~=~ [m_{2},m_{1},m_{2,4,2}; \pm \ha  m_{4} ]\nn\\ &&
 \chi^\pm_k ~=~ [m_{2,4},m_{1},m_{2,2}; 0  ]\nn
    \eea

These GVMs are related in the following way:
  \eqn{multop43} \begin{matrix} \L^-_b & {m_2\a_{24,34} \atop \longrightarrow } &\L^-_c &
   {m_1\a_{14,34} \atop \longrightarrow } & \L^-_d &  {m_2\a_{14,24} \atop \longrightarrow } & \L^-_j &  {m_4\a_{24,3} \atop \longrightarrow }  &\L^-_k \cr
&&&&\cr \updownarrow  && \updownarrow && \updownarrow && \updownarrow && \parallel  \cr
 &&&&\cr   \Lp_b &{m_2\a_{14,14,23,3} \atop \longleftarrow } & \Lp_c & {m_1\a_{14,24,23,3} \atop \longleftarrow }
 & \Lp_d &  {m_2\a_{14,24,3,3} \atop \longleftarrow } & \Lp_j &  {m_4\a_{44} \atop \longleftarrow } &\Lp_k
\end{matrix}
\ee
Note that the Remark before \eqref{multopjk} is again valid for the pair $\chi^\pm_j$.

\bigskip

The reduced multiplets of type $M_4$ also contain 15 GVMs (ERs).  Those induced  by finite-dimensional representations
  of $\cm$~ are:
\eqnn{mult44} &&\chi^-_a ~=~ \{ m_1,m_2,m_{3}, 0\}= \chi^-_0 \\ %
     && \chi^-_e ~=~ \{ m_{13}, m_2, m_{3}, -2m_{23}\} = \chi^-_g\nn\\ &&
\chi^-_f ~=~ \{ m_{23},m_{12},m_{3},-2m_{13}\} = \chi^-_i \nn\\ &&
\chi^-_k ~=~ \{m_{23},m_{1},m_{23,2},-2m_{13} -2m_2   \} = \chi^-_h \nn\\ &&
\chi^-_j ~=~ \{m_{2},m_{1},m_{23,23},-m_{13,13,23,2}     \} = \chh_j \nn\\ &&
 \chh_k  ~=~ \{ m_{23},m_{1},m_{23,2},-m_{13,13,23,2}  \}= \chh_h  \nn\\ && 
\chh_f   ~=~ \{ m_{23},m_{12},m_{3},-m_{13,13,23,2}  \}= \chh_i \nn\\ && 
 \chh_e  ~=~ \{ m_{13},m_{2},m_{3},-m_{13,13,23,2}  \} = \chh_g  \nn\\ && 
 \chh_a ~=~ \{ m_{1},m_{2},m_{3},-m_{13,13,23,2}  \} = \chh_0  \nn  
         \eea
or
\eqnn{mult4c4}
&&\chi^\pm_a ~=~ [ m_1,m_2,m_{3}; \pm (m_{13,2} + m_3/2)   ] \\
&&
\chi^\pm_e ~=~ [ m_{13}, m_2, m_{3};  \pm (m_{12} + m_3/2)  \nn\\ &&
\chi^\pm_f ~=~ [ m_{23},m_{12},m_{3}; \pm (m_{2} + m_3/2)  ]\nn\\ &&
\chi^\pm_k ~=~ [m_{23},m_{1},m_{23,2}; \pm \ha  m_{3} ]\nn\\ &&
 \chi^\pm_j ~=~ [m_{2},m_{1},m_{23,23}; 0  ] \nn
     \eea

These GVMs are related in the following way:
  \eqn{multop44} \begin{matrix} \L^-_a & \phantom{\lra} & \L^-_e & {m_1\a_{14,34} \atop \longrightarrow }
  & \L^-_f &  {m_2\a_{14,24} \atop \longrightarrow } & \L^-_k &  {m_3\a_{14} \atop \longrightarrow } & \L^-_j  \cr
&&&& && && \cr \updownarrow  &&\updownarrow  && \updownarrow && \updownarrow &&\parallel   \cr
 &&&&&& && \cr  \Lp_a && \Lp_e & {m_1\a_{14,24,23,3} \atop \longleftarrow } &
 \Lp_f & {m_2\a_{14,24,3,3} \atop \longleftarrow } & \Lp_k &  {m_3\a_{24,3} \atop \longleftarrow } & \Lp_j
\end{matrix}
\ee
 Again the Remark before \eqref{multopjk} is  valid for the pair $\chi^\pm_k$.

\subsection{Further reduced multiplets}

We list only those multiplets that contain GVMs (ERs) that are induced  by finite-dimensional representations
  of $\cm$~:

\bigskip

Type M13:
\eqnn{mult413} &&
\chi^-_c ~=~ \{ m_{2},m_{2},m_{4},-2m_{2}-m_{4}  \} \nn\\ && 
  \chh_c =\{ m_{2},m_{2},m_{4},-4m_{2}-2m_{4}  \} \nn\\ && 
    \chi^\pm_c ~=~ [ m_{2},m_{2},m_{4}; \pm ( m_{2} + m_4/2) ]
      \eea
  The two ERs of the above doublet are related only by the KS operators.

\bigskip

  Type M24:
  \eqnn{mult424}
     &&  \chi^-_k ~=~ \{ m_{3},m_{1},m_{3},-2m_{1}-2m_{3}\} \\ && 
     \chh_k  ~=~ \{ m_{3},m_{1},m_{3},-2m_{1}-3m_{3}   \}   \nn\\ &&
\chi^\pm_k ~=~ [m_{3},m_{1},m_{3}; \pm \ha  m_{3} ]
\eea
 Again the Remark before \eqref{multopjk} is   valid for the pair $\chi^\pm_k$.

\bigskip

  Type M34: Finally, we have a singlet:
\eqn{mult434}
\chi^\pm_k ~=~ \{m_{2},m_{1},2m_{2},-2m_{1}-4m_{2}    \}=   [m_{2},m_{1},2m_{2};  0 ]
     \ee

\section{Concluding remarks}

Matters are arranged so that in every main multiplet only the ER with
signature ~$\chi_0^-$~ contains a finite-dimensional nonunitary
subrepresentation in  a finite-dimensional subspace ~$\ce$. The
latter corresponds to the finite-dimensional   irrep of ~$F''_4 $~ with
signature ~$\{ m_1\,, m_2\,, m_3\,, m_4 \}$. Thus, the main multiplets are in 1-to-1 correspondence
  with the finite-dimensional representations of $F''_4 $.

  The subspace ~$\ce$~ is annihilated by the
operator ~$G^+\,$,\ and is the image of the operator ~$G^-\,$. The
subspace ~$\ce$~ is annihilated also by the \ido{} ~$\cd_{m_4\a_4}$~ acting from
~$\chi^-_0$~ to ~$\chi^-_a\,$.
 When all ~$m_i=1$~ then ~$\dim\,\ce = 1$, and in that case
~$\ce$~ is also the trivial one-dimensional UIR of the whole algebra
~$\cg$. Furthermore in that case the conformal weight is zero:
~$d=\sevha+c_{\vert_{m_i=1}}=0$.

In the conjugate ER ~$\chi_0^+$~ there is a unitary discrete series
subrepresentation in  an infinite-dimen\-sional subspace  $\cd_0$. It
is annihilated by the operator  $G^- $,\ and is in the image of the
operator  $G^+ $ acting from ~$\chi_0^-$ and in the image of the
\ido\ ~$\cd^{m_4}_{\a_{14,23,3}}$~ acting from ~$\chi_a^+$.

Two more occurrences of discrete series are in the  infinite-dimen\-sional subspaces  ~$\cd_a$,~$\cd_b$ of the
ERs  ~$\chi_a^+$,  ~$\chi_b^+$, resp. As above they are annihilated by the operator  $G^- $,\ and are in the images of the
operator  $G^+ $ acting from ~$\chi_a^-$, ~$\chi_b^-$, resp. Furthermore the subspace  ~$\cd_a$~ is in the image
of the operator ~$\cd^{m_3}_{\a_{14,23}}$~ acting from ~$\chi_b^+$ and is annihilated by the
\ido\ ~$\cd^{m_4}_{\a_{14,23,3}}$.
 Furthermore the subspace  ~$\cd_b$~ is in the image
of the operator ~$\cd^{m_2}_{\a_{14,14,23,3}}$~ acting from ~$\chi_c^+$ and is annihilated by the
\ido\ ~$\cd^{m_3}_{\a_{14,23}}$.

\bigskip

 After the present paper the only  split rank one case that is not treated yet is $Sp(N,1)$\cite{Dobf}.

\bigskip\bigskip

\nt {\bf Acknowledgments.}~~
  The author has received partial support from  Bulgarian NSF Grant DN-18/1.


 \np

 \voffset -4cm

\begin{figure}[]
\begin{center}
\includegraphics[width=1.4\hsize]{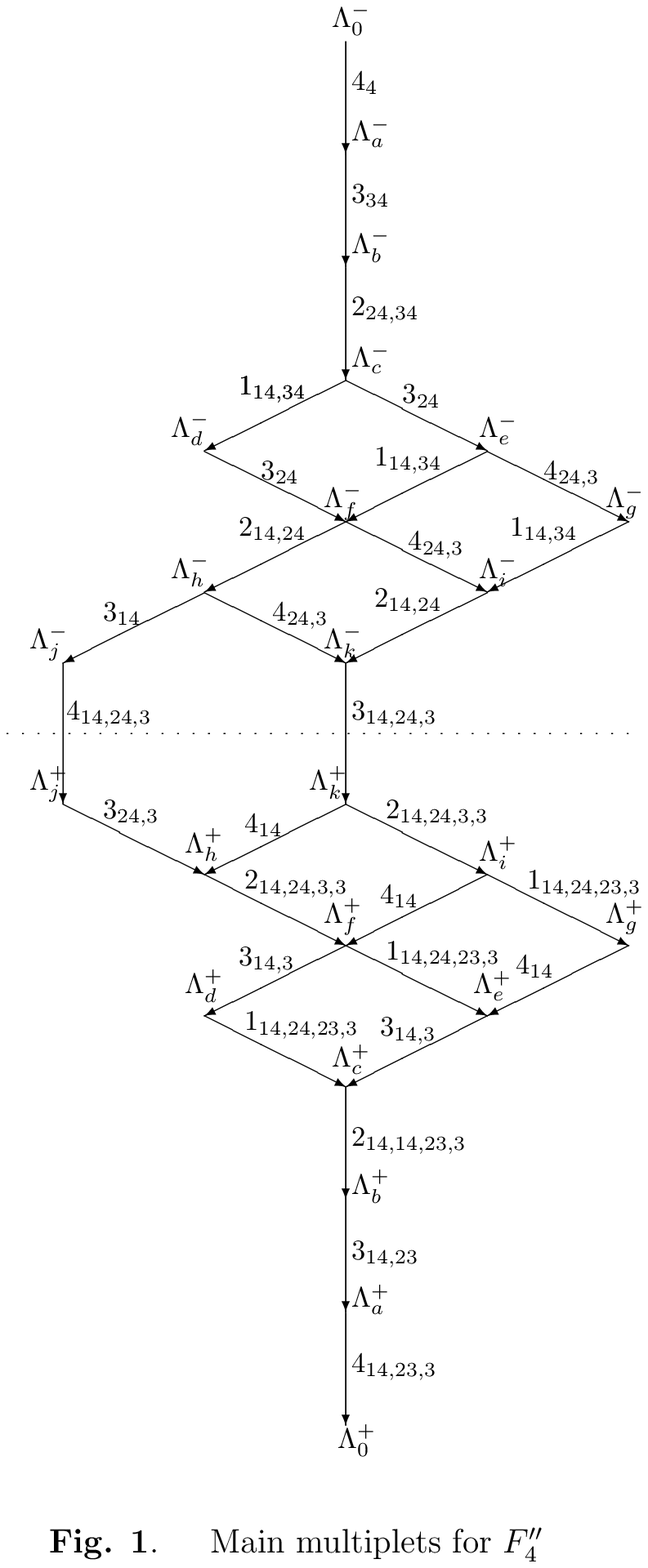} ~~~
\end{center}
\end{figure}



\end{document}